\documentclass{article}
\def\be{\begin{equation}}
\def\ee{\end{equation}}
\newcommand{\kk}[2]{\frac{#1}{#2}}

\newcommand{\ff}[1]{{\mbox{\boldmath $#1$}}}
\def\a{\alpha}
\def\b{\beta}

\def\={\approx}
\def\x{\ff{x}}
\def\vcode#1#2#3{\begin{figure}
\begin{center}
\begin{minipage}[c]{#1\textwidth}
\hrule \vspace{2pt}  %
{#2}  \vspace{5pt} \hrule
\end{minipage}
\caption{#3}
\end{center}   \end{figure}} 

\begin{document}

\title{Two-Stage Eagle Strategy with Differential Evolution}

\author{Xin-She Yang$^1$ and Suash Deb$^2$ \\ \\
1) Mathematics and Scientific Computing, \\
National Physical Laboratory, Teddington TW11 0LW, UK.
\and
2) Department of Computer Science \& Engineering, \\
C. V. Raman College of Engineering, \\
Bidyanagar, Mahura, Janla,
Bhubaneswar 752054, INDIA. }

\date{}

\maketitle

\begin{abstract}
Efficiency of an optimization process is largely determined by the search algorithm
and its fundamental characteristics. In a given optimization, a single type of
algorithm is used in most applications.  In this paper, we will investigate the
Eagle Strategy recently developed for global optimization, which uses
a two-stage strategy by combing two different algorithms to improve the overall
search efficiency. We will discuss this strategy with differential evolution and then
evaluate their performance by solving real-world optimization problems such as
pressure vessel and speed reducer design. Results suggest that we can reduce the computing
effort by a factor of up to 10 in many applications. \\ \\

\end{abstract}

\noindent Reference to this paper should be made as follows: \\ \\
Yang, X. S. and Deb, S., (2012).
`Two-Stage Eagle Strategy with Differential Evolution', \\
{\it Int. J. Bio-Inspired Computation}, Vol. 4, No. 1, pp.1--5.

\maketitle

\section{Introduction}

Metaheuristic optimization and computational modelling have become popular in
engineering design and industrial applications. The essence of such paradigm
is the efficient numerical methods and search algorithms. It is no exaggeration to
say that how numerical algorithms perform will largely determine the
performance and usefulness of modelling and optimization tools (Baeck et al.,1997;
Yang, 2010).

Among all optimization algorithms,  metaheuristic algorithms are becoming
powerful for solving tough nonlinear
optimization problems (Kennedy and Eberhart, 1995; Price et al., 2005; Yang, 2008; Cui and Cai, 2009).
The aim of developing modern metaheuristic algorithms is to enable the capability
of carrying out global search,  and good examples of nature-inspired
metaheuristics are particle swarm optimisation (PSO) (Kennedy and Eberhart, 1995)
 and Cuckoo Search (Yang and Deb, 2010a).
Most metaheuristic algorithms have relatively high efficiency in terms of finding
global optimality.

The efficiency of metaheuristic algorithms can be attributed to the
fact that they are designed to imitate the best features in nature,
especially the selection of the fittest in biological systems
which have evolved by natural selection over millions of years.
In real-world applications, most data have noise or associated randomness to a certain degree,
some modifications to these algorithms are often required,
in combination with some form of averaging or reformulation of the problem.
There exist some algorithms for stochastic optimization, and
the Eagle Strategy (ES), develop by Yang and Deb, is one of such algorithms for dealing with stochastic
optimization (Yang and Deb, 2010b).

In this paper, we will investigate the Eagle Strategy further by hybridizing it
with differential evolution (Storn, 1996; Storn and Price, 1997; Price et al., 2005). We  first validate the
ES by some multimodal test functions and then apply it to real-world optimization
problems. Case studies include pressure vessel design and gearbox speed reducer design.
We will discuss the results and point out directions for further research.

\section{Eagle Strategy}

Eagle strategy developed by Xin-She Yang and Suash Deb (Yang and Deb, 2010b)
is a two-stage method
for optimization. It uses a combination of crude global search and intensive local
search employing different algorithms to suit different purposes. In essence,
the strategy first explores the search space globally using a L\'evy
flight random walk, if it finds a promising solution, then an intensive
local search is employed using a more efficient local optimizer such as hill-climbing and downhill
simplex method. Then, the two-stage process starts again with
new global exploration followed by a local search in a new region.

The advantage of such a combination is to use a balanced tradeoff between
global search which is often slow and a fast local search. Some tradeoff and balance
are important. Another advantage of this method is that we can use any algorithms we like
at different stages of the search or even at different stages of iterations.
This makes it easy to combine the advantages of various algorithms so as
to produce better results.

\vcode{0.9}{Objective functions $f_1(\x), ..., f_N(\x)$ \\
Initialization and random initial guess $\x^{t=0}$ \\
{\bf while} (stop criterion) \\
Global exploration by randomization \\
Evaluate the objectives and find a promising solution \\
Intensive local search around a promising solution \\
\indent \qquad via an efficient local optimizer  \\
\indent \qquad {\bf if} (a better solution is found) \\
\indent \qquad \quad Update the current best \\
\indent \qquad {\bf end} \\
Update $t=t+1$ \\
{\bf end} \\
Post-process the results and visualization.}
{Pseudo code of the eagle strategy. \label{eagle-fig-50}  }

It is worth pointing that this is a methodology or strategy, not an algorithm.
In fact, we can use different algorithms at different stages and at different
time of the iterations.  The algorithm used for the global exploration should
have enough randomness so as to explore the search space diversely and
effectively. This process is typically slow initially, and should speed up
as the system converges (or no better solutions can be found after a certain
number of iterations). On the other hand, the algorithm used for the intensive
local exploitation should be an efficient local optimizer. The idea is to reach the
local optimality as quickly as possible, with the minimal number of function evaluations.
This stage should be fast and efficient.

\section{Differential Evolution}

Differential evolution (DE) was developed by R. Storn and K. Price by their
nominal papers  in 1996 and 1997 (Storn, 1996; Storn and Price, 1997).
It is a vector-based evolutionary
algorithm, and can be considered as a further development to genetic algorithms.
It is a stochastic search algorithm with self-organizing tendency
and does not use the information of derivatives. Thus, it is a population-based,
derivative-free method. Another advantage of differential evolution over
genetic algorithms is that DE treats solutions as real-number strings, thus
no encoding and decoding is needed.

As in genetic algorithms, design parameters in a $d$-dimensional search space
are represented as
vectors, and various genetic operators are operated over their bits of strings.
However, unlikely genetic algorithms, differential evolution carries out
operations over each component (or each dimension of the solution).
Almost everything is done in terms of vectors. For example,  in genetic
algorithms, mutation is carried out at one site or multiple sites of
a chromosome, while in differential evolution, a difference vector of two
randomly-chosen population vectors is used to perturb
an existing vector. Such vectorized mutation can be viewed as a
self-organizing search,
directed towards an optimality. This kind of perturbation is carried out over each
population vector, and thus can be expected to be more efficient. Similarly,
crossover is also a vector-based component-wise exchange of chromosomes
or vector segments.

For a $d$-dimensional optimization problem with $d$ parameters,
a population of $n$ solution vectors are initially generated,
we have $\x_i$ where $i=1,2,...,n$.
For each solution $\x_i$ at any generation $t$, we use the conventional
notation as
\be \x_i^t=(x_{1,i}^t, x_{2,i}^t, ..., x_{d,i}^t), \ee
which consists of $d$-components in the $d$-dimensional space.
This vector can be considered as the chromosomes or genomes.

Differential evolution consists of three main steps: mutation,
crossover and selection.

Mutation is carried out by the mutation scheme.
For each vector $\x_i$ at any time or generation $t$,
we first randomly choose three distinct
vectors $\x_p$, $\x_q$ and $\x_r$ at $t$,
and then generate a so-called donor vector
by the mutation scheme
\be \ff{v}_i^{t+1} =\x_p^t + F (\x_q^t-\x_r^t), \label{DE-F-equ-50} \ee
where $F \in [0,2]$ is a parameter, often referred to as
the differential weight. This requires that the minimum number
of population size is $n \ge 4$. In principle, $F \in [0,2]$,
but in practice, a scheme with $F \in [0,1]$ is more efficient and stable.
The perturbation $\ff{\delta}=F (\x_q-\x_r)$ to the vector $\x_p$
is used to generate a donor vector $\ff{v}_i$, and such perturbation
is directed and self-organized.

The crossover is controlled by a crossover probability $C_r \in [0,1]$
and actual crossover can be carried out in two ways: binomial and exponential.
The binomial scheme performs crossover on each of the $d$  components
or variables/parameters. By generating a uniformly distributed random number
$r_i \in [0,1]$, the $j$th component of $\ff{v}_i$ is manipulated as
\be \ff{u}_{j,i}^{t+1} =\ff{v}_{j,i} \quad \textrm{if } \; r_i \le C_r \ee
otherwise it remains unchanged.
This way, each component can be decided randomly whether to exchange with donor
vector or not.

Selection is essentially the same as that used in genetic algorithms. It is to select
the most fittest, and for minimization problem, the minimum objective value.

Most studies have focused on the choice of $F$, $C_r$ and $n$ as well as
the modification of (\ref{DE-F-equ-50}).  In fact, when generating mutation vectors,
we can use many different ways of formulating (\ref{DE-F-equ-50}), and this
leads to various schemes with the naming convention: DE/x/y/z where
x is the mutation scheme (rand or best), y is the number of difference vectors,
and z is the crossover scheme (binomial or exponential).
The basic DE/Rand/1/Bin scheme is given in (\ref{DE-F-equ-50}).
For a detailed review on different schemes, please refer to Price et al. (2005).

\section{ES with DE}

As ES is a two-stage strategy, we can use different algorithms at different stage.
The large-scale coarse search stage can use randomization via L\'evy flights
In the context of metaheuristics,
the so-called L\'evy distribution is a distribution of the
sum of $N$ identically and independently distribution random variables (Gutowski, 2001;
Mantegna, 1994; Pavlyukevich, 2007).

This  distribution is defined by a Fourier transform in  the following form
\be F_N(k) = \exp [-N |k|^{\b}].  \ee
The inverse to get the actual distribution $L(s)$ is not straightforward, as the integral
\be L(s) = \kk{1}{\pi} \int_0^{\infty} \cos (\tau s) e^{-\a \; \tau^{\b}} d\tau, \quad (0 < \b \le 2), \ee
does not have analytical forms, except for a few special cases. Here $L(s)$ is called
the L\'evy distribution with an index $\b$. For most applications, we can set $\a=1$ for
simplicity. Two special cases are $\b=1$ and $\b=2$. When $\b=1$, the above integral becomes the
Cauchy distribution.
When $\b=2$, it becomes the normal distribution. In this case, L\'evy flights become the
standard Brownian motion.

For the second stage, we can use differential evolution as the intensive local search.
We know DE is a global search algorithm, it can easily be tuned to do efficient
local search by limiting new solutions locally around the most promising region.
Such a combination may produce better results than those
by using pure DE only, as we will demonstrate this later. Obviously,
the balance of local search (intensification) and global search (diversification) is
very important, and so is the balance of the first stage and second stage in the ES.

\section{Validation}

Using our improved ES with DE, we can first validate it against some test functions which
are highly nonlinear and multimodal.

There are many test functions, here we have chosen the following 5 functions as a subset
for our validation.

Ackley's function
\[ f(\x) = -20 \exp\Big[-\kk{1}{5} \sqrt{\kk{1}{d}
\sum_{i=1}^d x_i^2} \; \Big]  \]
\be -\exp\Big[\kk{1}{d} \sum_{i=1}^d \cos ( 2 \pi x_i)  \Big]  + 20 +e, \ee
where $d=1,2,...$, and $-32.768 \le x_i \le 32.768$ for $i=1,2,...,d$. This function
has the global minimum $f_*=0$ at $\x_*=(0,0,...,0)$. \\

The simplest of De Jong's functions is the so-called sphere function
\be f(\x) =\sum_{i=1}^d x_i^2, \quad
-5.12 \le x_i \le 5.12, \ee
whose global minimum is obviously $f_*=0$ at $(0,0,...,0)$. This function is unimodal and convex.

Rosenbrock's function \be f(\x) = \sum_{i=1}^{d-1} \Big[ (x_i-1)^2 + 100 (x_{i+1}-x_i^2)^2 \Big], \ee
whose global minimum $f_*=0$ occurs at $\x_*=(1,1,...,1)$ in the domain
$-5 \le x_i \le 5$ where $i=1,2,...,d$. In the 2D case, it is often written as
\be f(x,y)=(x-1)^2 + 100 (y-x^2)^2, \ee
which is often referred to as the banana function.

Schwefel's function \be f(\x) = - \sum_{i=1}^d x_i \sin \Big(\sqrt{|x_i|} \Big),
\;\; -500 \le x_i \le 500, \ee
whose global minimum $f_*\=-418.9829 n$ occurs at $x_i=420.9687$ where $i=1,2,...,d$.

Shubert's function \[ f(\x)=\Big[ \sum_{i=1}^K i \cos \Big( i + (i+1) x \Big) \Big] \]
\be \cdot \Big[ \sum_{i=1}^K i \cos \Big(i + (i+1) y \Big) \Big], \ee
which has multiple global minima $f_* \= -186.7309$ for $K=5$ in the search domain $-10 \le x,y \le 10$.

Table I summarizes the results of our simulations, where $9.7\%$ corresponds to
the ratio of the number of function evaluations
in ES to the number of function evaluations in DE. That is the computational effort in ES
is only about $9.7\%$ of that using pure DE. As we can see that ES with DE is significantly
better than pure DE.

\begin{table}[ht]
\caption{Ratios of computational time}
\centering
\begin{tabular}{lllll}
\hline \hline
Functions & ES/DE  \\
\hline
Ackley ($d=8$)  & 24.9\%   \\
De Jong ($d=16$) & 9.7\% \\
Rosenbrock ($d=8$) & 20.2\% \\
Schwefel ($d=8$)  & 15.5\% \\
Shubert &  19.7\% \\
\hline
\end{tabular}
\end{table}

\section{Design Benchmarks}
Now we then use the ES with DE to solve some real-world case studies
including pressure vessel and speed reducer problems.

\subsection{Pressure Vessel Design}

Pressure vessels are literally everywhere such as champagne bottles and gas tanks.
For a given volume and working pressure, the basic aim of designing
a cylindrical vessel is to minimize the total cost. Typically, the design variables are
the thickness $d_1$ of the head, the thickness $d_2$ of the body, the inner radius
$r$, and the length $L$ of the cylindrical section (Coello, 2000; Cagnina et al., 2008).
This is a well-known test problem for optimization  and it can be written as
\[ \textrm{minimize } f(\x) = 0.6224 d_1 r L + 1.7781 d_2 r^2 \]
\be + 3.1661 d_1^2 L + 19.84 d_1^2 r, \ee
subject to the following constraints
\be
\begin{array}{lll}
 g_1(\x) = -d_1 + 0.0193 r \le 0 \\
 g_2(\x) = -d_2 + 0.00954 r \le 0 \\
 g_3(\x) = - \pi r^2 L -\frac{4 \pi}{3} r^3 + 1296000 \le 0 \\
 g_4(\x) =L -240 \le 0.
\end{array}
\ee

The simple bounds are  \be 0.0625 \le d_1, d_2 \le 99 \times 0.0625, \ee
and \be 10.0 \le r, \quad L \le 200.0. \ee

\begin{table}[ht]
\caption{Comparison of number of function evaluations }
\centering
\begin{tabular}{lllll}
\hline \hline
Case study & Pure DE & ES & ES/DE  \\
\hline
Pressure vessel  & 15000 & 2625 & 17.7\%   \\
Speed reducer & 22500 & 3352 & 14.9\% \\

\hline
\end{tabular}
\end{table}

Recently, Cagnina et al. (2008) used an efficient particle swarm optimiser
to solve this problem and they found the best solution
$f_* \approx 6059.714$
at \be \x_* \approx (0.8125, \; 0.4375, \; 42.0984, \; 176.6366). \ee
This means the lowest price is about $\$6059.71$.

Using ES, we obtained the same results, but we used significantly
fewer function evaluations, comparing using pure DE and other methods.
This again suggests ES is very efficient.

\subsection{Speed Reducer Design}

Another important benchmark is the design of a speed reducer
which is commonly used in many mechanisms such as
a gearbox (Golinski, 1973).
This problem involves the optimization of 7 variables, including
the face width, the number of teeth,
the diameter of the shaft and others. All variables are continuous within some limits,
except $x_3$ which only takes integer values.

\[ f(\x) = 0.7854 x_1 x_2^2 (3.3333 x_3^2+14.9334 x_3-43.0934) \]
\[ -1.508 x_1 (x_6^2+x_7^2)+7.4777 (x_6^3+x_7^3)  \]
\be +0.7854 (x_4 x_6^2+x_5 x_7^2) \ee

\be g_1(\x) = \frac{27}{x_1 x_2^2 x_3}-1 \le 0, \ee
\be g_2(\x) = \frac{397.5}{x_1 x_2^2 x_3^2}-1 \le 0 \ee

\be g_3(\x)=\frac{1.93 x_4^3}{x_2 x_3 x_6^4} - 1 \le 0, \ee
\be g_4(\x)=\frac{1.93 x_5^3}{x_2 x_3 x_7^4} - 1 \le 0 \ee

\be g_5(\x) =\frac{1.0}{110 x_6^3} \sqrt{\Big(\frac{745.0 x_4}{x_2 x_3} \Big)^2+16.9 \times 10^6} -1 \le 0 \ee

\be g_6(\x) =\frac{1.0}{85 x_7^3} \sqrt{\Big(\frac{745.0 x_5}{x_2 x_3} \Big)^2+157.5 \times 10^6} -1 \le 0 \ee

\be g_7(\x) = \frac{x_2 x_3}{40} -1 \le 0 \ee
\be g_8(\x) = \frac{5 x_2}{x_1} -1 \le 0 \ee

\be g_9(\x)=\frac{x_1}{12 x_2} -1 \le 0 \ee
\be g_{10}(\x) =\frac{1.5 x_6+1.9}{x_4}-1 \le 0 \ee
\be g_{11}(\x) =\frac{1.1 x_7+1.9}{x_5}-1 \le 0 \ee

where the simple bounds are $ 2.6 \le x_1 \le 3.6$, $0.7 \le x_2 \le 0.8$,
$17 \le x_3 \le 28$, $7.3 \le x_4 \le 8.3$, $7.8 \le x_5 \le 8.4$, $2.9 \le x_6 \le 3.9$,
and $5.0 \le x_7 \le 5.5$.

In one of latest studies, Cagnina et al. (2008) obtained
the following solution
\be \x_*=(3.5, 0.7, 17, 7.3, 7.8, 3.350214,5.286683) \ee
with $f_{\min}=2996.348165.$

Using our ES, we have obtained the new best
\be \x_*=(3.5, 0.7, 17, 7.3, 7.8, 3.34336449,5.285351) \ee
with the best objective
$f_{\min}=2993.7495888$.
We can see that ES not only provides better solutions but also finds
solutions more efficiently using fewer function evaluations.

\section{Discussions}

Metaheuristic algorithms such as differential evolution and eagle strategy are
very efficient. We have shown that a proper combination of these two can produce
even better performance for solving nonlinear global optimization problems.
First, we have validated the ES with DE and compared their performance.
We then used them to solve real-world optimization problems including
pressure vessel and speed reducer design. Same or better results have
been obtained, but with significantly less computational effort.

Further studies can focus on the sensitivity studies of the parameters
used in ES and DE so as to identify optimal parameter ranges for
most applications. Combinations of ES with other algorithms may also
prove fruitful. Furthermore, convergence analysis can provide
even more profound insight into the working of these algorithms.


\end{document}